\documentclass{elsart3-1}

\usepackage{amssymb}
\usepackage{amsbsy,amsmath,amsfonts,amssymb,amscd}
\usepackage{latexsym,euscript,times}
\usepackage[english,francais]{babel}

\newcommand{\PP}{\mathbb{P}}

\newcommand{\E}{\mathcal E}
\newcommand{\F}{\mathcal F}

\newtheorem{e-proposition}[theorem]{Proposition}

\newtheorem{e-definition}[theorem]{Definition\rm}

\newtheorem{theoreme}{Th\'eor\`eme}[section]

\newtheorem{remarque}{\it Remarque}
\newtheorem{exemple}{\it Exemple\/}
\renewcommand{\theequation}{\arabic{equation}}
\def\og{\leavevmode\raise.3ex\hbox{$\scriptscriptstyle\langle\!\langle$~}}
\def\fg{\leavevmode\raise.3ex\hbox{~$\!\scriptscriptstyle\,\rangle\!\rangle$}}

\begin{document}

\begin{frontmatter}




 
%
\selectlanguage{francais}
\title{Exemples de surfaces canoniques de ${\PP}^{6}$ et de solides de Calabi-Yau de ${\PP}^{7}$
}




\author[authorlabel1]{Marie-Am\' elie Bertin}
\ead{lbertin@math.unizh.ch, marie.bertin@parisfree.com}

\address[authorlabel1]{Universit\" at Z\" urich, Institut f\" ur Mathematik, Wintherturerstrasse, 190 , CH-8057 Z\" urich, SUISSE}

\begin{abstract}
Using  the global Gulliksen-Neg{\aa}rd complex, we build  
regular canonical surfaces of general type in ${\PP}^6$, 
Calabi-Yau $3$-folds in ${\PP}^{7}$ and Fano anticanonical $4$-folds,
all of degree $17$ and $20$. We also give their Hodge numbers and syzygies.

\vskip 0.5\baselineskip

\selectlanguage{francais}
 
\noindent{\bf R\'esum\'e}
Nous construisons \` a l'aide du complexe de Gulliksen-Neg{\aa}rd global des exemples 
de surfaces lisses canoniques r\' eguli\` eres, de solides de Calabi-Yau de ${\PP}^{7}$ 
et de vari\' et\' es de Fano anticanoniques lisses de dimension $4$ , 
de degr\' e $17,20$. Nous donnons leurs nombres de Hodge et syzygies.   
\vskip 0.5\baselineskip
\noindent

\end{abstract}

\end{frontmatter}


\selectlanguage{francais}
\section{Introduction et rappels}
\label{}
Soit $X$ un sous-sch\' ema localement de Cohen-Macaulay de codimension $c$ d'une vari\' et\' e complexe, projective et lisse ${\PP}$.
Le sch\' ema $X$ est dit \emph{sous-canonique} s'il existe un fibr\' e en droites 
${\mathcal{N}}$ sur $X$ tel que le faisceau dualisant
${\omega}_{X}$ de $X$ soit isomorphe comme ${{\mathcal O}}_{X}$-module \` a ${\mathcal{N}}\otimes {{\mathcal O}}_{X}$. La question de la d\' etermination des
types de r\' esolution libre  possible des sch\' emas sous-canoniques de ${\PP}$ pour une codimension fix\' ee n'est r\' esolue compl\` etement que pour 
la codimension $2$ (\cite{S},\cite{H}) mais a cependant conduit \` a la d\' ecouverte de nombres de complexes de longueur $c$, exacts lorsqu'ils
d\' eterminent des sous-sch\' emas de la bonne codimension $c$. En retour, lorsque ${\PP}$ est un espace projectif ${\PP}^r$,
ces constructions permettent de construire des surfaces canoniques et des solides de Calabi-Yau (\cite{O},\cite{T} et \cite{R}).
Parmi les complexes r\' esolvant \' eventuellement des sous-sch\' emas de codimension $4$, le plus ancien  est celui de Gulliksen et Neg{\aa}rd \cite{GN},
il r\' esoud dans sa version locale l'id\' eal des mineurs sous-maximaux d'une matrice carr\'ee.

Soient ${\E}$ et ${\F}$ des fibr\' es vectoriels sur ${\PP}$ de m\^ eme rang $e\geq 3$. Consid\' erons $\phi$ un morphisme de fibr\' es vectoriels
${\E}\xrightarrow{}{\F}$. Nous noterons simplement  par ${{\mathcal O}}$ le faisceau structural de ${\PP}$ et par ${{\mathcal L}}$ le fibr\' e vectoriel en droites
${\bigwedge}^{e}{\E}\otimes {\bigwedge}^{e}{\F}$.
Rappelons qu'il est associ\' e \` a ces donn\' ees un complexe $\mathbb{G}_{\bullet}$ de fibr\' es vectoriels 
sur ${\PP}$ dont la localisation en tout point $x\in {\PP}$ est isomorphe \` a un complexe de Gulliksen-Neg{\aa}rd  \cite{La} (complexe scandinave).
Remarquons que Pragacz et Weyman \cite{PW} donnent une construction ind\' ependante 
du choix d'une base donn\' ee  du complexe de Gulliksen-Neg{\aa}rd. Elle se globalise imm\' ediatement  pour donner le complexe global.  Associons \` a $\phi$ le dual $s_{\phi}^{*}$ de la section naturelle  associ\' ee \` a ${\bigwedge}^{e-1}\phi$ et
d\' efinissons le sous-sch\' ema de Gulliksen- Neg{\aa}rd $X(\phi )$   par ${{\mathcal O}}_{X(\phi )}=coker s_{\phi}^{*}$. 
Le complexe de Gulliksen-Neg{\aa}rd global $\mathbb{G}_{\bullet}$ a la forme suivante
$$0\xrightarrow{}{{\mathcal L}}^{\otimes 2}\xrightarrow{} {\E}\otimes {\F}^{*} \otimes L \xrightarrow{} det({\E})\otimes {\bigwedge}_{1,e-1} {\F}^{*}\oplus {det {\F}^{*}}\otimes {\bigwedge}_{1,e-1} {\E}\xrightarrow{} {\bigwedge}^{e-1} {\E}\otimes {\bigwedge}^{e-1}{\F}^{*}\xrightarrow{s_{\phi}^{*}}{{\mathcal O}}$$
o\` u ${\bigwedge}_{1,e-1}$ d\' esigne le foncteur de Schur associ\' e \` a la partition adjointe de $(e-1,1)$.

\section{Sous-variet\' es de Gulliksen-Negard}

Le complexe $\mathbb{G}_{\bullet}$ \' etant localement isomorphe au complexe de 
Gulliksen-Neg{\aa}rd, on d\' eduit la propri\' et\' e suivante de \cite{GN} (th\' eor\` eme 4 ) et du lemme d'acyclicit\' e sur les anneaux locaux de Cohen-Macaulay 

\begin{theoreme} Soient  ${\PP}$ une vari\' et\' e projective lisse et ${\E}$ et ${\F}$ deux fibr\' es vectoriels sur ${\PP}$ de m\^ eme rang $e\geq 3$.
Choisissons  $\phi\in Hom({\E},{\F})$. 
\begin{enumerate}

\item Si $X(\phi)$ est de codimension $4$, alors ${\mathbb{G}}_{\bullet}$ est une r\' esolution de $X(\phi)$.

\item Si de plus $X$ est Gorenstein, $X(\phi)$ est sous-canonique avec ${\omega}_{X(\phi )}\simeq {\omega}_{{\PP}}\otimes {{\mathcal L}}^{-\otimes 2}|_{X}$.

\item Dans le cas particulier o\` u $X$ est l'espace projectif ${\PP}^{r}$ on obtient la formule suivante:
$$c_1 ({\omega}_{X(\phi)})=-2(c_1 ({\E})-c_{1}({\F}))-r-1,$$
si $c_1$ d\' esigne la premi\` ere classe de Chern d'un fibr\' e vectoriel.
\end{enumerate}
\end{theoreme}

L'assertion $3$ du th\' eor\` eme d\' ecoule de la propri\' et\' e suivante du complexe de Gulliksen-Neg{\aa}rd 
$$ \mathbb{G}_{\bullet}^{*}\simeq \mathbb{G}_{\bullet}\otimes {{\mathcal L}}^{\otimes 2}.$$  

Rappelons enfin dans le contexte des sch\' emas de Gulliksen-Neg{\aa}rd
les r\' esultats classiques sur les lieux de d\' eg\' enerescence que nous utiliserons par la suite.

\begin{theoreme}[Banic\v a; Fulton-Lazarsfeld] Soit ${\PP}$ une vari\' et\' e 
complexe projective lisse de dimension $r$ telle que $r\geq 5$.
Choisissons $\phi$ g\' en\' eriquement dans $Hom({\E},{\F})$. Supposons de plus que ${\E}^{*}\otimes{\F}$ est engendr\' e par ses sections.

\begin{enumerate}

\item  Le sous-sch\' ema $X(\phi )$ est r\' eduit et irr\' eductible; lorsqu'il n'est pas vide il est de codimension $4$. 

\item  Si  $r\leq 8$ la vari\' et\' e $X(\phi)$ est vide ou lisse.

\item  Si ${E}^{*}\otimes {\F}$ est ample, $X(\phi)$ est non vide et lisse si et seulement si $r\leq 8$. 

\end{enumerate}

\end{theoreme}

En appliquant les th\' eor\` emes de Bertini affines (\cite{jou}, th\' eor\` eme 6.3 p. 66) 
 \` a la contruction que C. Banic\v a emploie, on d\' eduit sans mal l'assertion 1 du th\' eor\` eme. 
Les assertions suivantes combinent les r\' esultats de Banic\v a et Fulton-Lazarsfeld dans le cas des 
vari\' et\' es de Gulliksen-Neg{\aa}rd.

\section{Exemples dans ${\PP}^{r}$}

Dans cette section ${\PP}={\PP}^{r}$ pour $5\leq r\leq 8$ et l'on choisira toujours $\phi$ 
g\' en\' eriquement dans $Hom({\E},{\F})$. Rappelons que le degr\' e $d$ des varit\' es construites, 
lorsqu'elles sont de codimension $4$, s'obtient grace \` a la formule de Porteous, qui dans ce cas s'\' ecrit 
$d=c_2({\F}/{\E})^{2}-c_{1}({\F}/{\E})c_3({\F}/{\E})$, o\` u $c_i({\F}/{\E})$ d\' esigne le coefficient du terme
 de degr\' e $i$ du d\' eveloppement en s\' erie formelle du quotient de  polyn\^{o}mes de Chern $c({\F})/c({\E})$.

\begin{exemple} (Vari\' et\' es de Del Pezzo) Si l'on choisit $e=3$, ${\E}={{\mathcal O}}^{3}$ et  ${\F}={{\mathcal O}}^{3}(1)$ l'on obtient des vari\' et\' es de Del Pezzo
lisses de dimension $n\leq 4$. 
La r\' esolution de Gulliksen-Neg{\aa}rd est alors de la forme suivante:
$$0\xrightarrow{}{{\mathcal O}}(-6)\xrightarrow{}{{\mathcal O}}^{9}(-4)\xrightarrow{}{{\mathcal O}}^{16}(-3)\xrightarrow{}{{\mathcal O}}^{9}(-2)\xrightarrow{} {{\mathcal O}}.$$
\end{exemple}

Les exemples que nous construisons maintenant sont tels que ${\omega}_{X(\phi )}={{\mathcal O}}_{X(\phi )}(3-n)$, o\` u $n$ d\' esigne
la dimension de $X(\phi)$.

\begin{exemple} Si l'on choisit $e=4$, ${\E}={{\mathcal O}}^{4}$ et ${\F}={{\mathcal O}}^{4}(1)$ on obtient des surfaces lisses de degr\' e $20$ et de dimension $n\leq 4$.
Le complexe de Gulliksen-Neg{\aa}rd prend alors la forme suivante:
\begin{equation}
0\xrightarrow{}{{\mathcal O}}(-8)\xrightarrow{}{{\mathcal O}}^{16}(-5)\xrightarrow{}{{\mathcal O}}^{30}(-4)\xrightarrow{} {{\mathcal O}}^{16}(-3)\xrightarrow{} {{\mathcal O}}.\end{equation} 

On obtient ainsi des surfaces canoniques r\' eguli\` eres de degr\' e $20$, de genre sectionnel $\pi =21$ et de polyn\^ ome
de Hilbert $10x^2-10x+8$. 
On obtient aussi des solides de Calabi-Yau de ${\PP}^{7}$ de degr\' e $20$. Leur polyn\^ ome de Hilbert est $(10/3)x^3+(14/3)x$. 
Leurs losanges de Hodge sont respectivement
$$\begin{matrix}
&&1&&\\
&0&&0&\\
7&&60&&7\\
&0&&0&\\
&&1&&
\end{matrix}
\qquad \qquad
\begin{matrix}
&&&1&&&\\
&&0&&0&&\\
&0&&h^{1,1}&&0&\\
1&&34&&34&&1\\
&0&&h^{1,1}&&0&\\
&&0&&0&&\\
&&&1&&&
\end{matrix}$$

\underline{Autre construction:} Si l'on choisit $e=r$, ${\E}={{\mathcal O}}^{n+1}(1)\oplus {{\mathcal O}}^{3}$ et ${\F}={\mathcal T}$, o\` u ${\mathcal{T}}$ est le fibr\' e tangent sur ${\PP}$.  On obtient encore des vari\' et\' es de degr\' e $20$. 
Apr\` es saturation de leurs \' equations, leur r\' esolution libre est de m\^ eme forme que (1).
\end{exemple} 

\bigskip

\begin{exemple}
Si l'on choisit $e=3$, ${\E}={{\mathcal O}}^{3}$ et ${\F}={{\mathcal O}}^{2}(1)\oplus
{{\mathcal O}}(2)$, on obtient des vari\' et\' es de degr\' e $17$. Le complexe de
Gulliksen-Neg{\aa}rd est alors
\begin{equation}
0\xrightarrow{} {{\mathcal O}}(-8)\xrightarrow{}{{\mathcal O}}^{3}(-6)\oplus {{\mathcal O}}^{6}(-5)\xrightarrow{} {{\mathcal O}}^{12}(-4)\oplus {{\mathcal O}}^{2}(-3)\oplus {{\mathcal O}}^{2}(-5)\xrightarrow{} 
{{\mathcal O}}^{6}(-3)\oplus {{\mathcal O}}^{3}(-2)\xrightarrow{} {{\mathcal O}}
\end{equation} 

Les surfaces canoniques ainsi obtenues ont genre sectionnel $\pi =18$ et polyn\^ ome de
Hilbert $(17/2) x^{2}-(17/2)x +8$ et les solides de Calabi-Yau obtenus ont pour
polyn\^ ome de Hilbert $(17/6)x^{3}+(31/6)x$. Leurs losanges de Hodge 
sont respectivement
$$\begin{matrix}
&&1&&\\
&0&&0&\\
7&&63&&7\\
&0&&0&\\
&&1&&
\end{matrix}
\qquad \qquad
\begin{matrix}
&&&1&&&\\
&&0&&0&&\\
&0&&h^{1,1}&&0&\\
1&&58&&58&&1\\
&0&&h^{1,1}&&0&\\
&&0&&0&&\\
&&&1&&&
\end{matrix}
$$
\underline{Autre construction:} Si l'on choisit $e=r$, ${\E}={{\mathcal O}}^{n+2}(1)\oplus {{\mathcal O}}\oplus {{\mathcal O}}(-1)$ et ${\F}={\mathcal T}$ on obtient encore des vari\' et\' es de degr\' e $17$. Apr\` es saturation de leurs \' equations, leur r\' esolution libre est de la m\^ eme forme que (2).
\end{exemple}

\bigskip

\begin{remarque}
Ces exemples montrent l'existence de surfaces de type g\' en\' eral d'invariant $K^2 =20,17$, $p_g =7$ et $q=0$.
\end{remarque}

\bigskip

\begin{exemple} Si l'on choisit $e=r$, ${\E}={{\mathcal O}}^{n+3}(1)\oplus {{\mathcal O}}(-2)$ et ${\F}={\mathcal T}$, on obtient des vari\' et\' es de codimension $2$ dans
${\PP}^{r-2}$ obtenues par deux sections hyperplanes d'une intersection compl\` ete de deux hypersurfaces cubiques de ${\PP}^{r}$.

\bigskip

\end{exemple}

Pour effectuer les calculs n\' ec\' essaires nous avons utilis\' e le logiciel
de calcul formel Macaulay  2 \cite{M2}. Les logiciels de calculs formels ne peuvent manipuler que des modules libres;
aussi  dans le cas o\` u ${\F}={\mathcal T}$ utilise-t-on la repr\' esentation des morphismes entre puissances 
ext\' erieures du fibr\' e tangent \cite{DE} qui se d\' eduit du complexe de Koszul, pour calculer $s_{\phi}^{*}$.
Tous les calculs ont \' et\' e faits sur le corps fini ${\mathbb F}_{101}$.


\section*{Remerciements}
Je remercie les Professeurs C.Okonek et M.Brodmann pour leur soutien, et tout particuli\` erement le Prof. C. Okonek 
pour m'avoir sugger\' e cette \' etude. Cette note a \' et\' e influenc\' ee par des discussion avec le Professeur F.-O. Schreyer
sur un sujet voisin. Je l'en remercie  chaleureusement.

\end{document}